\documentclass[10pt]{amsart}
\usepackage{amsfonts,amssymb,amscd,amsmath,enumerate,verbatim,calc}
\usepackage{color,soul}
\usepackage{mathtools}
\usepackage{graphicx}
\everymath{\displaystyle}
\usepackage{fancyhdr}
\usepackage{tikz}
\usetikzlibrary{shapes.geometric, arrows}
\tikzstyle{arrow} = [thick,->,>=stealth]
\tikzstyle{process} = [rectangle, minimum width=4cm, minimum height=2cm, text centered, text width=2.5cm, draw=green, fill=green!10]
\theoremstyle{plain}

\theoremstyle{definition}

\begin{document}
\title{Q*-normal spaces }
\author{Hamant Kumar}
\email{hamant1984@gmail.com}
\address{Department of Mathematics, V. A. Government Degree College, Atrauli, Aligarh, 202280, U. P. (India)}

\author{Neeraj Kumar Tomar}
\email{neer8393@gmail.com}
\address{Department of Applied Mathematics, Gautam Buddha University, Greater Noida, Uttar Pradesh 201312, India}

\subjclass[IMS]{Primary 54D15, 54D10, 54A05, 54C08}
\date{\today}
	
\begin{abstract}
In this paper, using $Q^*$-closed sets, we introduce a new version of normality called, $Q^*$-normality which is a weak form of normality. Further utilizing $Q^*g$-closed sets, we obtain some characterizations of $Q^*$-normal and normal spaces and also obtain some preservation theorems for $Q^*$-normal spaces.
\end{abstract}
\maketitle

\section*{Introduction} 
In 1968, the notion of quasi normal space was introduced by Zaitsev\cite{zaitsev1968certain}. In 1970, Levine\cite{levine1970generalized} initiated the study of so called generalized closed (briefly $g$-closed) sets in order to extend many of the most important properties of closed sets to a large family. In 1973, Singal and Singal\cite{singal1973mildly} introduced the concept of mildly normal spaces and obtained their characterizations. In 1986, Munshi\cite{munshi1986separation} introduced and studied notions of $g$-normal and $g$-regular spaces. In 1990, Lal and Rahman\cite{lal1990note} have further studied notions of quasi normal and mildly normal spaces. In 2000, Dontchev and Noiri\cite{dontchev2000quasi} introduced the notion of $\pi g$-closed sets and using $\pi g$-closed sets,  obtained a new characterization of quasi normal spaces. In 2010, M. Murugalingam and N. Lalitha\cite{murugalingam2010q} introduced the concept of $Q^*$-closed sets and obtained some properties. In 2015, P. Padma and S. Udayakumar\cite{padma2015q} introduced the concept of $Q^*g$-closed sets and obtained some basic properties of $Q^*g$-closed sets.  In 2018, H. Kumar\cite{kumar2018} introduced and studied various forms of normal spaces in topological spaces in his Ph. D. Thesis. In 2024, H. Kumar et. al.\cite{kumar2024} introduced $Fg$-closed sets and obtained some properties of F-normal spaces in topological spaces in terms of $Fg$-closed sets. Recently In 2024, H. Kumar\cite{kumarQg2024} further studied $Q^*g$-closed sets and obtained a result.
\section{Preliminaries and Notations} 
Throughout this paper, spaces $(X,\tau)$, $(Y,\sigma)$, and $(Z,\rho)$ (or simply $X$, $Y$ and $Z$) always mean topological spaces on which no separation axioms are assumed unless explicitly stated. Let $A$ be a subset of a space $X$. The closure of $A$ and interior of A are denoted by $cl(A)$ and $int(A)$  respectively. A subset $A$ of a topological space $(X,\tau)$ is said to be \textbf{regular open (resp. regular closed)}  if $A= int(ci(A))$ (resp. $A = cl(int(A))$). 

\subsection{Definition} A subset $A$ of a space $(X,\tau)$ is said to be $Q^*$-closed\cite{murugalingam2010q} if $int(A)$ $=$ $\phi$ and $A$ is closed. The complement of a $Q^*$-closed set is said to be $Q^*$-open.

\subsection{Definition} A subset A of a space $(X,\tau)$ is said to be

$(i)$ \textbf{generalized closed (briefly g-closed)}\cite{levine1970generalized} if $cl(A)\subset U$ whenever $A\subset U$ and $U\in \tau$.

$(ii)$ \textbf{$Q^*g$-closed}\cite{kumarQg2024} if $cl(A)\subset U$ whenever $A\subset U$ and $U$ is $Q^*$-open in $X$.

\vspace{.1cm}
The complement of $g$-closed (resp. $Q^*g$-closed set) is called $g$-open (resp. $Q^*g$-open).

\vspace{.1cm}

The family of all $Q^*$-closed (resp. $Q^*$-open, $Q^*g$-closed, $Q^*g$-open) sets of a space $X$ is denoted by $Q^*$-$C(X)$ (resp. $Q^*$-$O(X)$, $Q^*g$-$C(X)$, $Q^*g$-$O(X)$).

\subsection{Remark.} We have the following implications for the properties of subsets:

\vspace{4.0mm}

\begin{flushleft}
\hspace{3cm}regular closed \\ \vspace{.5cm} \hspace{4cm}$\Downarrow$ \\ \vspace{.5cm}  
$Q^*$-closed\hspace{.5cm} $\Longrightarrow$ \hspace{1cm}  closed\hspace{.5cm} $\Longrightarrow$\hspace{.6cm} $g$-closed\hspace{.5cm} $\Longrightarrow$\hspace{.6cm} $Q^*g$-closed \\ \vspace{.2cm}
\end{flushleft}

Where none of the implications is reversible as can be seen from the following examples:
\subsection{Example} Let $X=\{a,b,c\}$ and $\tau = \{\phi,\{a\},\{b\},\{a,b\}, X\}$, Then

$(i)$ closed set are : $\phi$, $X$, $\{c\}$ $\{a,c\}$, $\{b,c\}$.

$(ii)$ $Q^*$-closed set are : $\phi$, $\{c\}$.

$(iii)$ $g$-closed set are : $\phi$, $X$, $\{c\}$, $\{a,c\}$, $\{b,c\}$.

$(iv)$ $Q^*g$-closed set are : $\phi$, $X$, $\{b\}$ $\{c\}$, $\{a,c\}$, $\{b,c\}$.

$(v)$ regular closed set are : $\phi$, $X$, $\{c\}$, $\{a,c\}$, $\{b,c\}$.

In the above example, every $g$-closed set is $Q^*g$-closed but converse is not true. Then the set $A = \{b\}$ is $Q^*g$-closed but not $g$-closed.

\subsection{Example} In $R$ with usual metric, finite sets are $Q^*$-closed but not regular closed. $[0, 1]$ is regular closed but not $Q^*$-closed. Hence regular closed and $Q^*$-closed sets are independent of each other.

\subsection{Theorem} For $Q^*g$-closed sets of a space $X$, the following properties hold:

\begin{flushleft}
    $(i)$ Every finite union of $Q^*g$-closed sets is always a $Q^*g$-closed.
\end{flushleft}

\begin{flushleft}
    $(ii)$ Every finite intersection of $Q^*g$-closed sets is always a $Q^*g$-closed.
\end{flushleft}

\subsection{Lemma} If $A$ be a subset of $X$, then

$(i)$ $Q^*$-$cl(X-A)$ = $X-Q^*$-$int(A)$.

$(ii)$ $Q^*$-$int(X-A)$ = $X-Q^*$-$cl(A)$.

\subsection{Theorem} A subset A of a space $X$ is $Q^*g$-open iff $F\subset int(A)$ whenever $F$ is $Q^*$-closed and $F\subset A$.

\textbf{Proof}. Let $F$ be $Q^*$-closed set such that $F\subset A$. Since $X–A$ is $Q^*g$-closed and $X-A \subset X-F$ where $F\subset int(A)$.

\textbf{Conversely}. Let $F\subset int(A)$ where $F$ is $Q^*$-closed and $F\subset A$. Since $F\subset A$ and $X–F$ is $Q^*$-open, $cl(X–A)$ $=$ $X –int(A)\subset X–F$. Therefore $A$ is $Q^*g$-open.

\section{$Q^*$-Normal Spaces}

\subsection{Definition} A space $X$ is said to be $Q^*$-normal if for every pair of disjoint  $Q^*$-closed subsets $A$, $B$ of $X$, there exist disjoint open sets $U$, $V$ of $X$ such that $A\subset U$ and  $B\subset V$.

\subsection{Definition} A space $X$ is said to be $g$-normal\cite{munshi1986separation} if for every pair of disjoint  $g$-closed subsets $A$, $B$ of $X$, there exist disjoint open sets $U$, $V$ of $X$ such that $A\subset U$ and  $B\subset V$.

\subsection{Definition} A space $X$ is said to be \textbf{quasi-normal}\cite{lal1990note} if for every pair of disjoint $\pi$-closed subsets $A$, $B$ of $X$, there exist disjoint open sets $U$, $V$ of $X$ such that $A\subset U$ and $B\subset V$.

\subsection{Definition} A space $X$ is said to be \textbf{midly-normal}\cite{singal1973mildly} if for every pair of disjoint regular closed subsets $A$, $B$ of $X$, there exist disjoint open sets $U$, $V$ of $X$ such that $A\subset U$ and $B\subset V$.

\textbf{By the definitions stated above, we have the following diagram:}
\vspace{5.0mm}

\begin{flushleft}
  \hspace{1cm} $g$-normal \hspace{.5cm} $\Longrightarrow$\hspace{.6cm} normal\hspace{.5cm} $\Longrightarrow$\hspace{.6cm}  $Q^*$-normal \\ 
\vspace{.4cm} \hspace{4.9cm} $\Downarrow$ \\ \vspace{.4cm} \hspace{4cm}  quasi-normal\hspace{.5cm} $\Longrightarrow$\hspace{.6cm} midly-normal 
\end{flushleft}
\vspace{5.0mm}
Where none of the implications is reversible as can be seen from the following examples:

\subsection{Example} Let $X=\{a,b,c,d\}$ and $\tau = \{\phi,\{a\}, \{c\},\{a,c\},\{a,b,d\},\{b,c,d\}, X\}$, Then the space $X$ is normal as well as $Q^*$ -normal. 

\subsection{Example} In $R$ with usual metric, finite sets are \textbf{$Q^*$-closed but not regular closed}. $[0, 1]$ is \textbf{regular closed but not $Q^*$-closed}. Hence regular closed and $Q^*$-closed sets are independent of each other.\\

From the above example, we can say that \textbf{mildly normal and $Q^*$-normal spaces} are independent of each other.

\subsection{Theorem} For a topological space $X$ the following properties are equivalent:\\
$(1)$ $X$ is $Q^*$-normal.\\
$(2)$ For every pair of $Q^*$-open subsets $U$ and $V$ of $X$ whose union is $X$, there exist closed subsets $G$ and $H$ of $X$ such that $A\subset U$, $B\subset V$ and $A\cup B = X$.\\
$(3)$ For any $Q^*$-closed set $A$ and every $Q^*$-open set $B$ in $X$ such that $A\subset B$, there exists an open subset $U$ of $X$ such that$A\subset U\subset Q^*$-$cl(U)\subset B$.\\
$(4)$ For every pair of disjoint $Q^*$-closed subsets $A$ and $B$ of $X$, there exist open subsets $U$ and $V$ of $X$ such that $A\subset U$, $B\subset V$ and $cl(U)\cap cl(V)$ = $\phi$.

\vspace{1.7mm}
\begin{flushleft}
\textbf{Proof}:  $(1)\Longrightarrow (2)$,  $(2)\Longrightarrow (3)$,  $(3)\Longrightarrow (4)$ and  $(4)\Longrightarrow (1)$

$(1)\Longrightarrow (2)$: Let $U$ and $V$ be any $Q^*$-open subsets of a $Q^*$-normal space $X$ such that $U\cup V = X$. Then, $X-U$ and $X-V$ are disjoint $Q^*$-closed subsets of $X$. By $Q^*$-normality of $X$, there exist disjoint open subsets $U_1$ and $V_1$ of $X$ such that $X-U\subset U_1$ and  $X-V\subset V_1$. Let $G = X-U_1$ and $H = X-V_1$.  Then, $G$ and $H$ are closed subsets of $X$ such that $G\subset U$, $H\subset V$ and $G\cup H = X$.
\end{flushleft}
\vspace{2.9mm}
\begin{flushleft}

$(2)\Longrightarrow(3)$: Let $A$ be a $Q^*$-closed and $B$ is a $Q^*$-open subsets of $X$ such that $A\subset B$. Then, $X-A$ and $B$ are $Q^*$-open subsets of $X$ such that $(X-A)\subset B = X$. Then, by part (2), there exist closed subsets $G$ and $H$ of $X$ such that $G\subset (X-A)$, $H\subset B$ and $G\cup H = X$. Then, $A\subset (X-G)$, $(X-B)\subset (X-H)$ and $(X-G)\cap (X-H) =\phi$. Let $U = X–G$ and $V = X-H$. Then $U$ and $V$ are disjoint open sets such that $A\subset U\subset X-V\subset B$. Since $X-V$ is closed, then we have $cl(U)\subset (X-V)$. Thus, $A\subset U\subset cl(U)\subset B$.

\end{flushleft}
\vspace{2.9mm}
\begin{flushleft}
$(3)\Longrightarrow(4)$: Let $A$ and $B$ be any disjoint $Q^*$-closed subset of $X$. Then $A\subset X-B$, where $X–B$ is $Q^*$-open. By the part (3), there exists an open subset $U$ of $X$ such that $A\subset U\subset cl(U)\subset X-B$. Let $V = X-cl(U)$. Then, $V$ is an open subset of $X$. Thus, we obtain $A\subset U$, $B\subset V$ and $cl(U)\cap cl(V)=\phi$.

$(4)\Longrightarrow(1)$: It is obvious.
\end{flushleft}

\subsection{Theorem} For a topological space $X$, the following properties are equivalent:\\
$(1)$ $X$ is $Q^*$-normal.\\
$(2)$ for any disjoint $H$, $K\in Q^*$-$C(X)$, there exist disjoint $Q^*g$-open sets $U$, $V$ such that $H\subset U$ and $K\subset V$.\\
$(3)$ for any $H\in Q^*$-$C(X)$ and any $V\in Q^*$-$O(X)$ containing $H$, there exists a $Q^*g$-open set $U$ of $X$ such that $H\subset U\subset Q^*g$-$cl(U)\subset V$.\\
$(4)$ for any $H\in Q^*$-$C(X)$ and any $V\in Q^*$-$O(X)$ containing $H$, there exists an open set $U$ of $X$ such that $H\subset U\subset Q^*g$-$cl(U)\subset V$.\\
$(5)$ for any disjoint $H$, $K\in Q^*$-$C(X)$, there exist disjoint open sets $U$, $V$ such that $H\subset U$ and $K\subset V$.

\vspace{1.7mm}

\textbf{Proof}:  $(1)\Longrightarrow (2)$: Since every open set is $Q^*g$-open, the proof is obvious. 

\begin{flushleft}
$(2)\Longrightarrow (3)$: Let $H\in Q^*$-$C(X)$ and $V$ be any $Q^*$-open set containing $H$. Then $H$, $X–V\in Q^*$-$C(X)$ and $H\cap (X–V)$ = $\phi$. By $(2)$, there exist $Q^*g$-open sets $U$, $G$ such that $H\subset U$, $X–V\subset G$ and $U\cap G = \phi$. Therefore, we have $H\subset U\subset (X–G) \subset V$. Since $U$ is $Q^*g$-open and $X–G$ is $Q^*g$-closed, we obtain $H\subset U\subset Q^*g$-$cl(U)\subset (X–G)\subset V$. 
\end{flushleft}

\begin{flushleft}
    
$(3)\Longrightarrow (4)$: Let $H\in Q^*$-$C(X)$ and $H\subset V\in Q^*$-$O(X)$. By $(3)$, there exist a $Q^*g$-open set $U_0$ of $X$ such that $H\subset U_0\subset Q^*g$-$cl(U_0)\subset V$. Since $Q^*g$-$cl(U_0)$ is $Q^*g$-closed and $V\in Q^*$-$O(X)$, $cl(Q^*g$-$cl(U_0))\subset V$. Put $int(U_0)$ = $U$, then $U$ is open and $H\subset U\subset cl(U)\subset V$. 
\end{flushleft}

\begin{flushleft}
    $(4)\Longrightarrow (5)$ Let $H$, $K$ be disjoint $Q^*$-closed sets of $X$. Then $H\subset (X–K)\in Q^*$-$O(X)$ and by $(4)$ there exists an open set $U_0$ such that $H\subset U_0\subset cl(U_0)\subset (X–K)$. Therefore, $V_0$ =$(X$ – $cl(U_0))$ is an open set such that $H\subset U_0$, $K\subset V_0$ and $U_0\cap V_0 = \phi$ . Moreover, put $U$ = $int(cl(U_0))$ and 
    
    $V$ = $int (cl(V_0))$, then $U$, $V$ are regular open sets such that $H\subset U$, $K\subset V$ and $U\cap V = \phi$ 

    $(5)\Longrightarrow(1)$: It is obvious.

\end{flushleft}

\subsection{Theorem} For a topological space $X$, the following properties are equivalent:\\
$(1)$ $X$ is normal.\\
$(2)$ for any disjoint closed sets $A$ and $B$, there exist disjoint $Q^*g$-open sets $U$ and $V$ such that $A\subset U$ and $B\subset V$.\\
$(3)$ for any closed set $A$ and any open set $V$ containing $A$, there exists a $Q^*g$-open set $U$ of $X$ such that $A\subset U\subset cl(U)\subset V$. 

\vspace{1.7mm}

\textbf{Proof}:  $(1)\Longrightarrow (2)$: Since every open set is $Q^*g$-open, the proof is obvious. 

\begin{flushleft}
$(2)\Longrightarrow (3)$: Let $A$ be a closed set and $V$ be any open set containing $A$. Then $A$ and $(X–V)$ are disjoint closed sets. There exist disjoint $Q^*g$-open sets $U$ and $W$ such that $A\subset U$ and $(X–V)\subset W$. Since $X–V$ is closed, we have $(X–V)\subset int(W)$ and $U\cap int(W) =\phi $ . Therefore, we obtain $cl(U)\cap int(W) =\phi $ and hence $A\subset U\subset cl(U)\subset (X–int(W))\subset V$. 
\end{flushleft}

\vspace{1.7mm}
\begin{flushleft}
$(3)\Longrightarrow (1)$ Let $A$, $B$ be disjoint closed sets of $X$. Then $A\subset (X–B)$ and $(X–B)$ is open. By $(3)$, there exists a $Q^*g$-open set $G$ of $X$ such that $A\subset G\subset cl(G)\subset (X–B)$. Since $A$ is closed, we have $A\subset int(G)$. Put $U = int(G)$ and $V = (X–cl(G))$. Then $U$ and $V$ are disjoint open sets of $X$ such that $A\subset U$ and $B\subset V$. Therefore, $X$ is normal. 
\end{flushleft}

\subsection{Proposition} Let $f : X\longrightarrow Y$ be a function, then:\\
$(i)$ The image of open subset under an open continuous function is open.\\
$(ii)$ The inverse image of $Q^*$-open (resp. $Q^*$-closed) subset under an open continuous function is $Q^*$-open (resp. $Q^*$-closed) subset.\\
$(iii)$ The image of closed subset under an open and a closed continuous surjective function is open.

\subsection{Theorem} The image of a $Q^*$-normal space under an open continuous injective function is a $Q^*$-normal.
\vspace{1.7mm}
\begin{flushleft}
    \textbf{Proof} Let $X$ be a $Q^*$-normal space and let $f : X\longrightarrow Y$ be an open continuous injective function. We need to show that $f(X)$ is a $Q^*$-normal. Let $A$ and $B$ be any two disjoint $Q^*$-closed sets in $f(X)$. Since the inverse image of a $Q^*$-closed set under an open continuous function is a $Q^*$-closed. Then, $f^{-1}(A)$ and $f^{-1}(B)$ are disjoint $Q^*$-closed sets in $X$. By $Q^*$-normality of $X$, there exist open subsets $U$ and $V$ of $X$ such that $f^{-1}(A)\subset U$, $f^{-1}(B)\subset V$ and $U\cap V =\phi$. Since $f$ is an open continuous injective function, we have $A\subset f(U)$, $B\subset f(V)$ and $f(U)\cap f(V) =\phi $. By \textbf{Proposition 2.10}, we obtain $f(U)$ and $f(V)$ are disjoint open sets in $f(X)$ such that $A\subset f(U)$ and $B\subset f(V)$. Hence $f(X)$ is $Q^*$-normal.
\end{flushleft}
\subsection{Theorem} For a space $X$, the following are equivalent:\\
$(1)$ $X$ is $Q^*$-normal.\\
$(2)$ For any disjoint $Q^*$-closed sets $H$ and $K$, there exist disjoint $g$-open sets $U$ and $V$ such that $H\subset U$ and $K\subset V$.\\
$(3)$ For any disjoint $Q^*$-closed sets $H$ and $K$, there exist disjoint $Q^*g$-open sets $U$ and $V$ such that $H\subset U$ and $K\subset V$.\\
$(4)$ For any $Q^*$-closed set $H$ and any $Q^*$-open set $V$ containing $H$, there exists a $g$-open set $U$ of $X$ such that $H\subset U\subset cl(U)\subset V$.\\
$(5)$ For any $Q^*$-closed set $H$ and any $Q^*$-open set $V$ containing $H$, there exists a $Q^*g$-open set $U$ of $X$ such that $H\subset U\subset cl(U)\subset V$.

\vspace{1.7mm}
\begin{flushleft}
    \textbf{Proof}.$(1)\Longrightarrow (2)$, $(2)\Longrightarrow (3)$, $(3)\Longrightarrow (4)$, $(4)\Longrightarrow (5)$ and $(5)\Longrightarrow (1)$.
\end{flushleft}

\vspace{1.7mm}
\begin{flushleft}
$(1)\Longrightarrow (2)$. Let $X$ be $Q^*$-normal space. Let $H$, $K$ be disjoint $Q^*$-closed sets of $X$. By assumption, there exist disjoint open sets $U$, $V$ such that $H\subset U$ and $K\subset V$. Since every open set is $g$-open, so $U$ and $V$ are $g$-open sets such that $H\subset U$ and $K\subset V$.

\end{flushleft}
\vspace{1.7mm}
\begin{flushleft}
    $(2)\Longrightarrow (3)$. Let $H$, $K$ be two disjoint $Q^*$-closed sets. By assumption, there exist disjoint $g$-open sets $U$ and $V$ such that $H\subset U$ and $K\subset V$. Since every $g$-open set is $Q^*g$-open, so $U$ and $V$ are $Q^*g$-open sets such that $H\subset U$ and $K\subset V$.
\end{flushleft}
\vspace{1.7mm}

\begin{flushleft}
    $(3)\Longrightarrow (4)$. Let $H$ be any $Q^*$-closed set and $V$ be any $Q^*$-open set containing $H$. By assumption, there exist disjoint $Q^*g$-open sets $U$ and $W$ such that $H\subset U$ and $X\subset V\subset W$. By \textbf{Theorem 1.8}, we get $X\subset V\subset int(W)$ and $cl(U)\cap int(W) =\phi$. Hence $H\subset U\subset cl(U)\subset X\subset int(W)\subset V$.
\end{flushleft}
\vspace{1.7mm}
\begin{flushleft}
    $(4)\Longrightarrow (5)$. Let $H$ be any $Q^*$-closed set and $V$ be any $Q^*$-open set containing $H$. By assumption, there exist $g$-open set $U$ of $X$ such that $H\subset U\subset cl(U)\subset V$. Since, every $g$-open set is $Q^*g$-open, there exists $Q^*g$-open set $U$ of $X$ such that $H \subset U\subset cl(U)\subset V$.
\end{flushleft}
\vspace{1.7mm}
\begin{flushleft}
    $(5)\Longrightarrow (1)$. Let $H$, $K$ be any two disjoint $Q^*$-closed sets of $X$. Then $H\subset X-K$ and $X-K$ is $Q^*$-open. By assumption, there exists $Q^*g$-open set $G$ of $X$ such that $H\subset G\subset cl(G)\subset X-K$. Put $U = int(G)$, $V = X-cl(G)$. Then $U$ and $V$ are disjoint open sets of $X$ such that $H\subset U$ and $K\subset V$.
\end{flushleft}

\section{Some Related Functions and $Q^*$-normal Spaces}

\subsection{Definition} A function $f : X\longrightarrow Y$  is said to be:\\
$(1)$ \textbf{almost $Q^*g$-continuous} if for any regular open set $V$ of $Y$, $f^{-1}(V)\in Q^*g$-$O(X$;\\ 
$(2)$ \textbf{almost $Q^*g$-closed} if for any regular closed set $F$ of $X$, $f(F)\in Q^*g$-$C(Y)$.

\subsection{Definition} A function $f : X\longrightarrow Y$  is said to be:\\
$(1)$ \textbf{$Q^*$-irresolute (resp. $Q^*$-continuous)} if for any $Q^*$-open (resp. open) set $V$ of $Y$, $f^{-1}(V)$ is $Q^*$-open in $X$;\\ 
$(2)$ \textbf{pre-$Q^*$-closed (resp. $Q^*$-closed)} if for any $Q^*$-closed (resp. closed) set $F$ of $X$, $f(F)$ is $Q^*$-closed in $Y$;

\subsection{Theorem} A function $f : X\longrightarrow Y$  is an almost $Q^*g$-closed surjection if and only if for each subset $S$ of $Y$ and each regular open set $U$ containing $f^{-1}(S)$, there exists a $Q^*g$-open set $V$ such that $S\subset V$ and $f^{-1}(V)\subset U$ . 

\vspace{1.7mm}
\begin{flushleft}
    \textbf{Proof}: \textbf{Necessity}. Suppose that $f$ is almost $Q^*g$-closed. Let $S$ be a subset of $Y$ and $U$ be a regular open set of $X$ containing $f^{-1}(S)$. Put $V =Y– f(X–U)$, then $V$ is a $Q^*g$-open set of $Y$ such that $S\subset V$ and $f^{-1}(V)\subset U$. 
\end{flushleft}

\vspace{1.7mm}
\begin{flushleft}
    \textbf{Sufficiency}. Let $F$ be any regular closed set of $X$. Then $f^{-1}(Y-f(F))\subset (X–F)$ and $X-F$ is regular open. There exists a $Q^*g$-open set $V$ of $Y$ such that $(Y-f(F))\subset V$ and $f^{-1}(V)\subset (X–F)$. Therefore, we have $f(F)\supset (Y–V)$ and $F\subset f^{-1}(Y–V)$. Hence, we obtain $f(F) = Y-V$ and $f(F)$ is $Q^*g$-closed in $Y$. This shows that $f$ is almost $Q^*g$-closed. 
\end{flushleft}

\subsection{Theorem} If $f : X\longrightarrow Y$  is an almost $Q^*g$-closed $Q^*$-irresolute (resp. $Q^*$-continuous) surjection and $X$ is $Q^*$-normal, then $Y$ is $Q^*$-normal (resp. normal).

\begin{flushleft}
   \textbf{Proof}. Let $A$ and $B$ be any disjoint $Q^*$-closed (resp. closed) sets of $Y$. Then $f^{-1}(A)$ and $f^{-1}(B)$ are disjoint $Q^*$-closed sets of $X$. Since $X$ is $Q^*$-normal, there exist disjoint open sets $U$ and $V$ of $X$ such that $f^{-1}(A)\subset U$ and $f^{-1}(B)\subset V$. Put $G = int(cl(U))$ and $H = int(cl(V))$, then $G$ and $H$ are disjoint regular open sets of $X$ such that $f^{-1}(A)\subset G$ and $f^{-1}(B)\subset H$. By \textbf{Theorem 3.3}, there exist $Q^*g$-open sets $K$ and $L$ of $Y$ such that $A\subset K$, $B\subset L$. $f^{-1}(K)\subset G$ and $f^{-1}(L)\subset H$. Since $G$ and $H$ are disjoint, so $K$ and $L$ are also disjoint. It follows from \textbf{Theorem 2.8 (resp. Theorem 2.9)} that $Y$ is $Q^*$-normal (resp. normal).
\end{flushleft}

\subsection{Theorem} If $f : X\longrightarrow Y$ is a continuous almost $Q^*g$-closed surjection and $X$ is a normal space, then $Y$ is normal.

\begin{flushleft}
    \textbf{Proof}. The proof is similar to that of \textbf{Theorem 3.4}. 

\end{flushleft}

\subsection{Theorem} : If $f : X\longrightarrow Y$ is an almost $Q^*g$-continuous pre-$Q^*$-closed (resp. $Q^*$-closed) injection and $Y$ is $Q^*$-normal, then $X$ is $Q^*$-normal (resp. normal).

\begin{flushleft}
    \textbf{Proof}.Let $H$ and $K$ be disjoint $Q^*$-closed (resp. closed) sets of $X$. Since $f$ is a pre-$Q^*$-closed (resp. $Q^*$-closed) injection, $f(H)$ and $f(K)$ are disjoint $Q^*$-closed sets of $Y$. Since $Y$ is $Q^*$-normal, there exist disjoint open sets $P$ and $J$ such that $f(H)\subset P$ and $f(K)\subset J$. Now, put $U = int(cl(P))$ and $V = int(cl(J))$, then $U$ and $V$ are disjoint regular open sets such that $f(H)\subset U$ and $f(K)\subset V$. Since $f$ is almost $Q^*g$-continuous, $f^{-1}(U)$ and $f^{-1}(V)$ are disjoint $Q^*g$-open sets such that $H\subset f^{-1}(U)$ and $K\subset f^{-1}(V)$. It follows from \textbf{Theorem 2.8 (resp. Theorem 2.9)} that $X$ is $Q^*$-normal (resp. normal). 
\end{flushleft}

\section{conflict of interest}
We certify that this work is original, has never been published before, and is n't being considered for publication anywhere else at this time. This publication is free from any conflicts of interest. As the Corresponding Author, I certify that each of the listed authors has read the paper and given their approval for submission.

\bibliographystyle{amsplain}
\bibliography{references}

\providecommand{\bysame}{\leavevmode\hbox to3em{\hrulefill}\thinspace}
\providecommand{\MR}{\relax\ifhmode\unskip\space\fi MR }
\providecommand{\MRhref}[2]{%
  \href{http://www.ams.org/mathscinet-getitem?mr=#1}{#2}
}
\providecommand{\href}[2]{#2}
\begin{thebibliography}{10}

\bibitem{dontchev2000quasi}
J~Dontchev and T~Noiri, \emph{Quasi-normal spaces and $\pi$ g-closed sets}, Acta Mathematica Hungarica \textbf{89} (2000), no.~3, 211--219.

\bibitem{kumar2018}
H.~Kumar, \emph{Some weaker forms of normal spaces in topological spaces}, Ph.d. thesis, C. C. S. University, Meerut, 2018.

\bibitem{kumar2024}
H.~Kumar, B.~S. Sharma, and A.~Kumar, \emph{Fg-closed sets and f-normal spaces}, International Journal of Creative Research Thoughts (IJCRT) \textbf{12} (2024), no.~4, I577--I582.

\bibitem{kumarQg2024}
Hamant Kumar, \emph{Q$^*g$-closed sets}, EPRA International Journal of Research and Development (IJRD) \textbf{9} (2024), no.~7, 189--190.

\bibitem{lal1990note}
Sunder Lal and MS~Rahman, \emph{A note of quasi-normal spaces}, Indian Journal of mathematics \textbf{32} (1990), no.~1, 87--94.

\bibitem{levine1970generalized}
Norman Levine, \emph{Generalized closed sets in topology}, Rendiconti del Circolo Matematico di Palermo \textbf{19} (1970), 89--96.

\bibitem{munshi1986separation}
BM~Munshi, \emph{Separation axioms}, Acta Ciencia Indica \textbf{12} (1986), no.~2, 140--145.

\bibitem{murugalingam2010q}
M~Murugalingam and N~Lalitha, \emph{Q* star sets}, Bulletin of Pure \& Applied Sciences-Mathematics and Statistics \textbf{29} (2010), no.~2, 369--376.

\bibitem{padma2015q}
P~Padma and S~UdayaKumar, \emph{Q* g-closed sets in topological spaces}, International Journal of advanced research in engineering and applied sciences \textbf{2} (2015), no.~1.

\bibitem{singal1973mildly}
MK~Singal and Asha~Rani Singal, \emph{Mildly normal spaces}, Kyungpook Mathematical Journal \textbf{13} (1973), no.~1, 27--31.

\bibitem{zaitsev1968certain}
V~Zaitsev, \emph{On certain classes of topological spaces and their bicompactifications}, Dokl. Akad. Nauk SSSR, vol. 178, 1968, pp.~778--779.

\end{thebibliography}

\end{document}